\documentclass[runningheads]{llncs}
\usepackage{graphicx}

\usepackage{amssymb,amsmath}
\usepackage{xcolor}
\usepackage{hyperref}
\usepackage[all]{xy}

\begin{document}
\title{Weighted homomorphisms on the \\
$p$-analog of the Fourier-Stieltjes algebra\\
 induced by piecewise affine maps}
\titlerunning{weighted homomorphisms of $B_p(G)$}
%

\author{Mohammad Ali Ahmadpoor\inst{1}\orcidID{0000-0001-6902-1916} \and Marzieh Shams Yousefi \inst{2}\orcidID{0000-0003-0426-708X}}
\authorrunning{M. A. Ahmadpoor \& M. Shams Yousefi}
%
\institute{{${}^{1}$ School of Mathematics and Statistics, Carleton University, Ottawa, Canada\\
${}^{2}$ Department of Pure Mathematics, Faculty of Mathematical Sciences, University of Guilan, Rasht, Iran\\
\email{mohammadaliahmadpoo@cmail.carleton.ca}\\
\email{m.shams@guilan.ac.ir}}}

\maketitle              
\begin{abstract}
In this paper, for $p\in(1,\infty)$ we study $p$-complete boundedness of weighted homomorphisms on the $p$-analog of the Fourier-Stieltjes algebras, $B_p(G)$, based on the $p$-operator space structure defined by the authors. Here,  for a locally compact group $G$, the space $B_p(G)$ stands for Runde's definition of the $p$-analog of the Fourier-Stieltjes algebra and the implemented $p$-operator space structure is come from the duality between $B_p(G)$ and the algebra of universal $p$-pseudofunctions, $UPF_p(G)$. It is established that the homomorphism $\Phi_\alpha:B_p(G)\to B_p(H)$, defined by $\Phi(u)=u\circ\alpha$ on $Y$ and zero otherwise, is $p$-completely contractive when the continuous and proper map $\alpha :Y\subseteq H\to G$ is affine, and it is $p$-completely bounded whenever $\alpha$ is piecewise affine map. Moreover, we assume that $Y$ belongs to the coset ring generated by open and amenable subgroups of $H$. To obtain the result, by utilizing the properties of $QSL_p$-spaces and representations on them, the relation between $B_p(G/N)$ and a closed subalgebra of $B_p(G)$ is shown, where $N$ is a closed normal subgroup of $G$. Additionally, $p$-complete boundedness of several well-known maps on such algebras are obtained.
\keywords{Completely bounded homomorphisms\and $p$-analog of the Fourier-Stieltjes algebras\and $QSL_p$-spaces\and Piecewise affine maps}        

\textbf{MSC2010:} Primary 46L07; Secondary 43A30, 47L10.
\end{abstract}

\section{Introduction}\label{subsection1.1}

Let $G$ be a locally compact group. The Fourier algebra, $A(G)$, and the Fourier-Stieltjes algebra, $B(G)$, on the locally compact group $G$, have been found by Eymard in 1964 \cite{EYMARD1964}. The general form of special type of maps on the Fourier and Fourier-Stieltjes algebras has been studied extensively. For example, when $G$ is an Abelian topological group, $A(G)$ is nothing except $L_1(\widehat{G})$, where $\widehat{G}$ is the Pontrjagin dual group of $G$, and $B(G)$ is isometrically isomorphic to $M(\widehat{G})$, the measure algebra. In this case, Cohen in \cite{COHEN1960-1} and \cite{COHEN1960-2} studied homomorphisms from $L_1(G)$ to $M(H)$, for Abelian groups $G$ and $H$, and gave the general form of these maps, as the weighted maps by a piecewise affine map on the underlying groups.\\
By \cite{BLECHER1992,EFFROSRUAN1991}, we know that $A(G)$ and $B(G)$ are operator spaces as the predual of a von Neumann algebra, and the dual of a $C^*$-algebra, respectively. Amini and the second author studied the relation between operator and order structure of these algebras via amenability \cite{SHAMSAMINISADY2010}. Ilie in \cite{ILIE2004} and \cite{ILIESPRONK2005} studied the completely bounded homomorphisms from the Fourier to the Fourier-Stieltjes algebras. It is shown that for a continuous piecewise affine map $\alpha:Y\subseteq H\rightarrow G$, the homomorphism $\Phi_\alpha:A(G)\rightarrow B(H)$, defined through
\begin{align*}
\Phi_\alpha u=\left\{
\begin{array}{ll}
u\circ\alpha& \text{on} \: Y\\
0& \text{o.w.}
\end{array}\right. ,\quad u\in A(G),
\end{align*}
is  completely bounded. Moreover, in the cases that $\alpha$ is an affine map and a homomorphism, the homomorphism $\Phi_\alpha$ is completely contractive and completely positive, respectively.

The Fig\`a-Talamanca-Herz algebras were introduced by Fig\`a-Talamanca for Abelian locally compact groups \cite{FIGATALAMANCA1965}, and it is generalized for any locally compact group by Herz \cite{HERZ1971}. For $p\in (1,\infty)$, coefficient functions of the left regular representation of a locally compact group $G$ on $L_{p}(G)$ make the Fig\`a-Talamanca-Herz algebra $A_p(G)$, and we have $A_2(G)=A(G)$. Therefore,  Fig\`a-Talamanca-Herz algebras can be considered as the $p$-analog of the Fourier algebras.

Daws in \cite{DAWS2010} introduced a $p$-operator space structure, with an extensive application to $A_p(G)$, which generalizes the natural and non-trivial operator space structure of $A(G)$.\\
Oztop and Spronk in \cite{OZTOPSPRONK2012}, and Ilie in \cite{ILIE2013} studied the $p$-completely bounded homomorphisms on the Fig\`a-Talamanca-Herz algebras, using this $p$-operator space structure. In \cite{ILIE2013} it is shown that the map $\Phi_\alpha :A_p(G)\rightarrow A_p(H)$, defined via
\begin{align*}
\Phi_\alpha u= \left\{
\begin{array}{ll}
u\circ \alpha & \text{on} \: Y\\
0 & \text{o.w}
\end{array}\right.,\quad u\in A_p(G),
\end{align*}
is a $p$-completely (bounded) contractive homomorphism for a continuous proper (piecewise) affine map $\alpha : Y\subseteq H\rightarrow G $ in the case that the locally compact group $H$ is amenable.

Runde in \cite{RUNDE2005} found a $p$-analog of the Fourier-Stieltjes algebras, $B_p(G)$. He used the theory of $ QSL_p $-spaces and representations on these spaces. Additionally, the $p$-operator space structure of $B_p(G)$ is fully described by the authors in \cite{AHSH2021-I} and it is shown that $B_p(G)$ is a $p$-operator space, as the dual space of the algebra of universal $p$-pseudofunctions $UPF_p(G)$. The second author of this paper studied the $p$-analog of the Fourier-Stieltjes algebras on the inverse semigroups in \cite{shams}.

In this paper, for a continuous proper piecewise affine map $\alpha : Y\subseteq H\rightarrow G $, we study the weighted map $\Phi_\alpha : B_p(G)\rightarrow B_p(H)$ which is defined by
\begin{align}\label{eq1}
\Phi_\alpha u= \left\{
\begin{array}{ll}
u\circ \alpha & \text{on} \: Y\\
0 & \text{o.w}
\end{array}\right.,\quad u\in B_p(G).
\end{align}
We will show that when $\alpha$ is an affine map, $\Phi_\alpha$ is a $p$-complete contraction, and in the case that $\alpha$ is a piecewise affine map, it is $p$-completely bounded homomorphism. At this aim, we put amenability assumption on open subgroups of $H$. Our approach to the concept of $p$-operator space structure on the $p$-analog of the Fourier-Stieltjes algebra, is the $p$-operator structure that can be implemented on this space from its predual.

The paper is organized as follows. In Section \ref{SECTIONPRELIMINARIES}, we give required  definitions and theorems about the $p$-analog of the Fourier-Stieltjes algebras and representations on $QSL_p$-spaces, and then some crucial and previously obtained properties of the algebra $B_p(G)$ will be listed. As a consequence, in Section \ref{SECTIONSPECIALMAPS} we will establish $p$-completely boundedness of some well-known operators on the $p$-analog of the Fourier-Stieltjes algebras (Theorem \ref{IMPORTANTMAPPING}). The obtained results will be applied in the final section, Section \ref{SECTIONINDUCEDHOMO}, to generalize Ilie's results on homomorphisms of the Fig\`a-Talamanca-Herz algebras in \cite{ILIE2013}.

\section{Preliminaries}\label{SECTIONPRELIMINARIES}
Herewith, we divide the prerequisites into three parts. First, we give general notions and theorems of representations and specific algebras on locally compact groups together with some useful tools to deal with our problem. Next, we focus on some properties of the $p$-analog of the Fourier-Stieltjes algebra and its $p$-operator space structure. Finally, we introduce some special maps regarding our main result. Throughout this paper, $G$ and $H$ are locally compact groups, and for $p\in (1,\infty)$, the number $p'$ is its complex conjugate, i.e. $1/p+1/p'=1$.

We commence by the essential description of $QSL_p$-spaces, and representations of groups on such spaces. For more information one can see \cite{RUNDE2005}.
\begin{definition}\label{def1}
A representation of a locally compact group $ G $ is a
pair $ (\pi , E) $, where $ E $ is a Banach space and $ \pi $ is a group homomorphism that maps each element $x\in G$ to an invertible isometric operator $\pi(x)$ on $E$. This homomorphism is continuous with respect to the given topology on $ G $
and the strong operator topology on $ \mathcal{B}(E) $.
\end{definition}

\begin{remark}
Each representation $ (\pi , E) $ can be lifted to a representation of the group algebra $L_1(G)$ on $ E $. Denoting this homomorphism with the same symbol $ \pi $, it is defined through
\begin{align}\label{11}
&\pi(f)=\int f(x)\pi (x)dx,\ f\in L_1(G),\\
&\langle \pi(f)\xi, \eta\rangle =\int f(x)\langle \pi (x)\xi, \eta\rangle dx,\quad \xi\in E,\:\eta\in E^*,\nonumber
\end{align}
where the integral \eqref{11} converges with respect to the strong operator topology.
\end{remark}

\begin{definition}
For two representations $ (\pi , E)$ and $(\rho , F) $ of the locally compact group $G$, we have the following terminologies.
\begin{enumerate}
\item
$ (\pi , E)$ and $ (\rho , F) $ are called equivalent, if there exists an invertible isometric map $T : E\rightarrow F$ for which the following diagram commutes for each $x\in G$,

\begin{displaymath}
\xymatrix{
E \ar[r]^{\pi(x)} & E \ar[d]^{T\quad .} \\
F \ar[r]_{\rho(x)} \ar[u]^{T^{-1}} & F }
\end{displaymath}

\item
 The representation $ (\pi , E)$ has a subrepresentation $ (\rho , F) $, if $F$ is a closed subspace of $E$, and for each $x\in G$, the operator $\rho(x)$ is the restriction of $\pi(x)$ to the subspace $F$.
\item
We say that $ (\pi , E)$ contains $ (\rho , F) $ and write $(\rho , F) \subseteq (\pi , E)$, if it is equivalent to a subrepresentation of $  (\pi , E) $.
\end{enumerate}
\end{definition}

\begin{definition}
\begin{enumerate}
\item A Banach space is called an $L_p$-space if it is of the form $L_p(X)$ for some measure
space $X$.
\item A Banach space is called a $ QSL_p $-space if it is isometrically isomorphic to a quotient of a subspace of an $ L_p $-space.
\end{enumerate}
\end{definition}
We denote by $ \text{Rep}_p(G) $ the collection of all (equivalence classes) of representations of $G $ on a  $ QSL_p $-space.

\begin{definition}
A representation of a Banach algebra $\mathcal{A}$ is a pair $ (\pi , E)$, where $E$ is a Banach space, and $\pi $ is a contractive algebra homomorphism from $\mathcal{A}$ to $\mathcal{B}(E)$. We call
$ (\pi , E)$ isometric if $\pi $ is an isometry and essential if the linear span of $\{\pi(a)\xi \: : \ a \in \mathcal{A},\: \xi\in
E\} $ is dense in $E$. 
\end{definition}

\begin{remark}\label{rem11}
For a locally compact group $G$ there exists a one-to-one correspondence between representations of $G$ and essential representations of $L_1(G)$.
\end{remark}

\begin{definition}

\begin{enumerate}
\item
If for a representation $ (\pi , E)\in \text{Rep}_p(G) $ there is an element $\xi_0\in E$ such that 
\begin{align*}
E=\overline{\{\pi (f)\xi_0\ : f\in L_1(G)\}{}}^{\|\cdot\|_E},
\end{align*}
then $(\pi,E)$ is said to be cyclic, and  the set of all such cyclic representations of the group $G$ on $QSL_p$-spaces is denoted by $\text{Cyc}_p(G)$.
\item
If a representation $ (\pi , E)\in \text{Rep}_p(G)$ contains (up to an isometry) each cyclic representation then it is called $p$-universal.
\end{enumerate}

\end{definition}

\begin{definition}
We say that the function $u:G\rightarrow\mathbb{C}$ is a coefficient function of a representation $(\pi,E)$, if there exist $\xi\in E$ and $\eta\in E^*$ such that
\begin{align*}
u(x)=\langle\pi(x)\xi,\eta\rangle,\quad x\in G.
\end{align*}
\end{definition}

The left regular representation of $G$ on $L_{p}(G)$, is denoted by  $ \lambda_{p}$ and defined as following
\begin{align*}
&\lambda_{p}: G\rightarrow \mathcal{B}(L_{p}(G)),\quad\lambda_{p}(x)\xi(y)=\xi(x^{-1}y),\quad \xi\in L_{p}(G),\:x,y\in G,
\end{align*}
and it obtains a description of the Fig\`a-Talamanca-Herz algebras.

\begin{definition}
\begin{enumerate}
\item
The set of all linear combinations of the coefficient functions of the representation $(\pi,E)\in\text{Rep}_p(G)$, namely elements of the form
\begin{align}\label{TTT1111}
u(x)=\sum_{n=1}^\infty\langle\pi(x)\xi_n,\eta_n\rangle,\quad x\in G,\; (\xi_n)_{n=1}^\infty\subseteq E,\quad (\eta_n)_{n=1}^\infty\subseteq E^*
\end{align}
where
\begin{align}\label{TTT222}
\sum_{n=1}^\infty\|\xi_n\|\|\eta_n\|<\infty,
\end{align}
together with the norm
\begin{align}\label{TTT333}
\|u\|=\inf\bigg\{\sum_{n=1}^\infty\|\xi_n\|\|\eta_n\|\ :\ u\;\text{as}\;\eqref{TTT1111}\;\text{with}\;\eqref{TTT222}\bigg\}
\end{align}
is denoted by $A_{p,\pi}$, and is called the $p$-analog of the $\pi$-Fourier space.
\item
In particular, if $(\pi,E)=(\lambda_p,L_p(G))$ then it is simply denoted by $A_p(G)$.
\item
More generally, if $(\pi,E)$ is a $p$-universal representation, then this space is denoted by $B_p(G)$.

\end{enumerate}
\end{definition}

\begin{remark}\label{REMARKRUNDERUNDE}
\begin{enumerate}
\item
By \cite[Theorem 4.7]{RUNDE2005}, the space $B_p(G)$ equipped with the norm defined as above, and pointwise operations is a commutative unital Banach algebra, and it is called the $p$-analog of the Fourier-Stieltjes algebra. Additionally, it has been shown that the space $A_p(G)$ is a Banach algebra which is called the Fig\`a-Talamanca-Herz algebra.

\item
The $p$-analog of the Fourier-Stieltjes algebra has been studied, for example in \cite{COWLING1979}, \cite{FORREST1994}, \cite{MIAO1996} and \cite{PIER1984}, as the multiplier algebra of the Fig\`a-Talamanca-Herz algebra. In this paper, we follow the construction of Runde in definition and notation (See \cite{RUNDE2005}) which we have swapped indexes $p$ and $p'$.

\item
By \cite[Corollary 5.3]{RUNDE2005}, if we denote the multiplier algebra of $A_p(G)$ by $\mathcal{M}(A_p(G))$, we have the following contractive embeddings
\begin{align*}
A_p(G)\subseteq B_p(G)\subseteq \mathcal{M}(A_p(G)).
\end{align*}

\item
One may express u in \eqref{TTT1111} as
\begin{align}\label{TTT444}
u(x)=\sum_{n=1}^\infty\langle\pi_n(x)\xi_n,\eta_n\rangle,\quad x\in G,\ \xi_n\in E_n,\ \eta_n\in E_n^*,
\end{align}
for which \eqref{TTT222} holds
and where $(\pi_n,E_n)\subseteq(\pi,E)$, and $(\pi_n,E_n)\in\text{Cyc}_p(G)$, for all $n\in \mathbb{N}$. In this case, the norm \eqref{TTT333} can be replaced by
\begin{align*}
\|u\|=\inf\bigg\{\sum_{n=1}^\infty\|\xi_n\|\|\eta_n\|\ :\ u\;\text{as}\;\eqref{TTT444}\;\text{with}\;\eqref{TTT222}\bigg\}
\end{align*}
\end{enumerate}
\end{remark}

\begin{definition}
\begin{enumerate}
\item
For a representation $(\pi,E)\in\text{Rep}_p(G)$, the closure of the space $\{\pi(f)\ :\ f\in L_1(G)\}$ with respect to the norm $\|\cdot\|_{\mathcal{B}(E)}$ is denoted by $PF_{p,\pi}(G)$ and it is called the algebra of $p$-pseudofunctions associated with $(\pi,E)$.
\item
If $(\pi,E)$ is a $p$-universal representation, then it is denoted by $UPF_p(G)$ and is called the algebra of
universal $ p $-pseudofunctions.
\item
If $(\pi,E)=(\lambda_p,L_p(G))$, then it is simply denoted by $PF_{p}(G)$.
\end{enumerate}
\end{definition}

Next lemma states that $ B_p(G) $  is a dual space.
\begin{lemma}{\cite[Lemma 6.5]{RUNDE2005}}\label{RUNDEDUALITY}
Let $ (\pi , E)\in \text{Rep}_{p}(G)$. Then, for each $  \phi\in PF_{p,\pi}(G)^*$, there is a unique $g\in B_p(G) $, with $\|g\|_{B_p(G)}\leq \|\phi\| $ such that
\begin{equation}\label{duality}
\langle\pi (f),\phi\rangle=\int_G f(x)g(x)dx,\qquad f\in L_1(G).
\end{equation}
Moreover, if $ (\pi , E) $ is $p$-universal, we have $\|g\|_{B_p(G)}= \|\phi\| $.
\end{lemma}

\begin{definition}
For a representation $(\pi, E)\in\text{Rep}_p(G)$ the dual space $PF_{p,\pi}(G)^*$ will be denoted by $B_{p,\pi}$, following the tradition initiated in \cite{ARSAC1976}, and it is called the $p$-analog of the $\pi$-Fourier-Stieltjes algebra.
\end{definition}

\begin{remark}
\begin{enumerate}
\item
The duality between $B_{p,\pi}$ and $PF_{p,\pi}(G)$ can be stated via relation below,
\begin{align*}
\langle \pi(f),u\rangle=\int_G u(x)f(x)dx,\quad f\in L_1(G),\ u\in B_{p,\pi}.
\end{align*}
\item
On top of that,
\begin{align*}
&\| u\|=\sup_{\|f\|_\pi\leq 1}|\langle \pi(f), u\rangle|=\sup_{\|f\|_\pi\leq 1}|\int_G u(x)f(x)dx|,\quad u\in B_{p,\pi},\\
&\| f\|_\pi=\sup_{\|u\|\leq 1}|\langle \pi(f), u\rangle|=\sup_{\|u\|\leq 1}|\int_G u(x)f(x)dx|,\quad f\in L_1(G).\\
\end{align*}
\item
Additionally, the fashionable notation for $PF_{p,\lambda_p}(G)^*$ is $PF_{p}(G)^*$, instead of $B_{p,\lambda_p}$. For the procedure of obtaining the space $B_{p,\pi}$ one would refer to \cite[Subsection 3.2]{AHSH2021-II}.

\end{enumerate}
\end{remark}
Now, we state a result from \cite{RUNDE2005} adapted to our notations.

\begin{theorem}\label{THEOREMRUNDE2005}
\begin{enumerate}
\item
For a representation $(\pi,E)\in\text{Rep}_p(G)$, the inclusion  $B_{p,\pi}\subseteq B_p(G)$ is contractive, and is an isometric isomorphism whenever $(\pi,E)$ is a $p$-universal representation.
\item
We have the following contractive inclusions
\begin{align*}
PF_{p}(G)^*=B_{p,\lambda_p}\subseteq B_p(G)\subseteq \mathcal{M}(A_p(G)),
\end{align*}
and all inclusions will become equalities in the case that $G$ is amenable.
\end{enumerate}
\begin{proof}
See \cite[Theorem 6.6 and Theorem 6.7]{RUNDE2005}).
\end{proof}
\end{theorem}

A concrete $p$-operator space is a closed subspace of $\mathcal{B}(E)$, for some $QSL_p$-space $E$. In this case for each $n\in\mathbb{N}$ one can define a norm $\|\cdot\|_n$ on $\mathbb{M}_n(X)=\mathbb{M}_n\otimes X$ by identifying $ \mathbb{M}_n(X) $ with a subspace of $ \mathcal{B}(l_p^n\otimes_p E $). So, we have the family of norms $\Big(\|\cdot\|_n\Big)_{n\in\mathbb{N}}$ satisfying:

\begin{enumerate}
\item[$\mathcal{D}_\infty:$] For $u\in \mathbb{M}_n(X)$ and $v\in \mathbb{M}_m(X)$, we have that $ \|u\oplus v\|_{n+m}=\max\{\|u\|_n,\|v\|_m\} $. Here $u\oplus v\in \mathbb{M}_{n+m}(X)$ has block representation
$\begin{pmatrix}
u & 0 \\
0 & v
\end{pmatrix}. $

\item[$\mathcal{M}_p:$] For every $u\in \mathbb{M}_m(X)$ and $\alpha\in \mathbb{M}_{n,m}$, $\beta\in \mathbb{M}_{m,n}$, we have that 
$$ \|\alpha u\beta\|_n\leq\|\alpha\|_{\mathcal{B}(l^m_p,l^n_p)}\|u\|_m\|\beta\|_{\mathcal{B}(l^n_p,l^m_p)}. $$

\end{enumerate}

Additionally, an abstract $p$-operator space is a Banach space $X$ equipped with the family of norms $(\|\cdot\|_n)_{n\in\mathbb{N}}$ defined by $\mathbb{M}_n(X)$ which satisfy two axioms above.

\begin{definition}
A linear operator $\Psi :X\rightarrow Y$ between two $p$-operator spaces is called $p$-completely bounded, if $\|\Psi\|_{\text{p-cb}}=\sup_{n\in\mathbb{N}}\|\Psi^{(n)}\|<\infty$, and $p$-completely contractive if $\|\Psi\|_{\text{p-cb}}=\sup_{n\in\mathbb{N}}\|\Psi^{(n)}\|\leq 1$, where $\Psi^{(n)} :\mathbb{M}_n(X)\rightarrow\mathbb{M}_n(Y)$ is defined in the natural way.
\end{definition}

\begin{lemma}\cite[Lemma 4.5]{DAWS2010}\label{DAWSLEMMA}
If $ \Psi:X\rightarrow Y $ is $p$-completely bounded map between two operator spaces $X$ and $Y$, then $\Psi^* :Y^*\rightarrow X^*$ is $p$-completely bounded, with $\|\Psi^*\|_{\text{p-cb}}\leq \|\Psi\|_{\text{p-cb}}$.
\end{lemma}

\begin{remark}
\begin{enumerate}
\item
It should be noticed that, converse of Lemma \ref{DAWSLEMMA} is not necessarily true, unless $X$ be a closed subspace of $\mathcal{B}(E)$, for some $L_p$-space $E$.
\item
In comparison with \cite{ILIE2013}, because of above explanations, a major difference in our work is that we need to study the predual maps of some crucial $p$-completely bounded maps (see Theorem \ref{IMPORTANTMAPPING}), instead of their dual maps, which is used in the classic theory.
\end{enumerate}
\end{remark}

Since the main purpose of this paper is to study $p$-completely boundedness of operators on the $p$-analog of the Fourier-Stieltjes algebras, we need to specify the $p$-operator structure on these algebras. We briefly explain about this structure in the following.

For a representation $(\pi,E)\in\text{Rep}_p(G)$, in \cite[Proposition 3.11]{AHSH2021-I}, it is indicated that the algebra $PF_{p,\pi}(G)$ is $p$-completely isomorph to a closed subspace of $\mathcal{B}(\mathcal{E})$, for some $QSL_p$-space $\mathcal{E}$. Indeed the space $\mathcal{E}$ has been constructed through $l_p$-direct sum of $QSL_p$-spaces associated with cyclic subrepresentations of $(\pi,E)$ in a particular process. As a consequence, it has been obtained that for a representation $(\pi,E)\in\text{Rep}_p(G)$, the algebra of $p$-pseudofunctions $PF_{p,\pi}(G)$ is a $p$-operator space \cite[Theorem 3.12]{AHSH2021-I}. Moreover, a $p$-operator space structure for the algebra $B_p(G)$ is derived from its predual $UPF_p(G)$. In fact, we have the following sequence of results from \cite{AHSH2021-I} in this regard.

\begin{theorem}{\cite[Theorem 4.5]{AHSH2021-I}}\label{TH45}
The algebra of universal $p$-pseudofunctions $UPF_p(G)$ has a non-trivial $p$-operator space structure.
\end{theorem}

\begin{theorem}\cite[Theorem 4.7]{AHSH2021-I}\label{TH47}
For $p\in (1,\infty)$, the Banach algebra $B_p(G)$ is a $p$-operator space.
\end{theorem}

Next proposition is the output of all materials mentioned before.

\begin{proposition}\cite[Proposition 4.8]{AHSH2021-I}\label{P48}
For a locally compact group $G$, and a complex number $p\in (1,\infty)$, the identification $B_p(G)=UPF_p(G)^*$ is $p$-completely isometric isomorphism.
\end{proposition}

In Section \ref{SECTIONINDUCEDHOMO}, we will study the homomorphisms on the $p$-analog of the Fourier-Stieltjes algebras induced by a continuous map $\alpha :Y\subseteq H\rightarrow G$, in the cases that $\alpha$ is homomorphism, affine and piecewise affine map, and $Y$ in the coset ring of $H$. So, we give some preliminaries here.\\
For a locally compact topological group $H$, let $\Omega_0(H)$ denote the ring of subsets which generated by open cosets of $H$. By \cite{ILIE2013} we have
\begin{align}\label{OPENCOSETRINGS}
\Omega_0(H)=\left\{Y\backslash\cup_{i=1}^nY_i \: :
\begin{array}{ll}
&Y\:\text{is an open coset of }\: H,\\
& Y_1,\ldots ,Y_n\:\text{open subcosets of infinite index in}\: Y\\
\end{array}\right\}.
\end{align}
Moreover, for a set $Y\subseteq H$, by $\text{Aff}(Y)$ we mean the smallest coset containing $Y$, and if $Y=Y_0\backslash\cup_{i=1}^nY_i\in\Omega_0(H)$, then $\text{Aff}(Y)=Y_0$.\\
Similarly, let us denote by  $\Omega_{\text{am-}0}(H)$ the ring of open cosets of open amenable subgroups of $H$, that is the ring of subsets of the form $Y\backslash\cup_{i=1}^nY_i$ where $Y$ is an open coset of an open amenable subgroup of $H$ and $Y_i$ is an open subcosets of infinite index in $Y$ (which has to be an open coset of an open amenable subgroup), for $i=1\ldots,n$.

\begin{definition}\label{DEFPIECEWISE}
Let $\alpha : Y\subseteq H\rightarrow G$ be a map.
\begin{enumerate}
\item
The map $\alpha$ is called an affine map on an open coset $Y$ of an open subgroup $H_0$, if
\begin{equation*}
\alpha(xy^{-1}z)=\alpha(x)\alpha(y)^{-1}\alpha(z),\qquad x,y,z\in Y,
\end{equation*}

\item
The map $\alpha$ is called a piecewise affine map if
\begin{enumerate}
\item
there are pairwise disjoint $ Y_i\in\Omega_0(H)$, for $ i=1,\ldots , n $, such that $Y=\cup_{i=1}^nY_i$,

\item there are affine maps $\displaystyle{\alpha_i : \text{Aff}(Y_i)\subseteq H\rightarrow G}$, for $ i=1,\ldots , n $, such that
\begin{equation*}
\alpha |_{Y_i}=\alpha_i |_{Y_i}.
\end{equation*}

\end{enumerate}
\end{enumerate}

\end{definition}

\begin{definition}
If $X$ and $Y$ are locally compact spaces, then a map $\alpha :Y\rightarrow X$ is called proper, if $\alpha^{-1}(K)$ is compact subset of $Y$, for every compact subset $K$ of $ X $.
\end{definition}

\begin{proposition}\cite[Proposition 4]{DUNKLRAMIREZ1971}\label{ramirez}
Let $\alpha : H\rightarrow G$ be a continuous group homomorphism. Then $\alpha$ is proper if and only if the bijective homomorphism $\tilde{\alpha}: H/{\ker\alpha}\rightarrow \alpha(H)=G_0$, is a topological group isomorphism, and $\ker\alpha$ is compact.
\end{proposition}

\begin{remark}\label{AFFFIINEREMMM}
\begin{enumerate}
\item
Proposition \ref{ramirez} implies that every continuous proper homomorphism is automatically a closed map. Therefore, $\alpha(H)$ is a closed subgroup of $G$. Additionally, $\ker\alpha$ is a compact normal subgroup of $H$.
\item
It is well-known that $\tilde{\alpha}$ is a group isomorphism, if and only if $\alpha$ is an open homomorphism into $\alpha(H)$, with the relative topology.  
\item\label{affine-remark}{\cite[Remark 2.2]{ILIE2004}}
If $Y=h_0H_0$ is an open coset of an open subgroup $H_0\subseteq H$, and $\alpha : Y\subseteq H\rightarrow G$ is an affine map, then there exists a group homomorphism $ \beta $ associated to $\alpha$ such that
\begin{align}\label{affine-homomorphism}
&\beta : H_0\subseteq H\rightarrow G,\quad\beta (h)= \alpha(h_0)^{-1}\alpha(h_0h),\quad h\in H_0.
\end{align}
\item \label{HOMAFFPROPER}
It is clear that, $\alpha$ is a proper affine map, if and only if $\beta$ is a proper homomorphism.
\item\cite[Lemma 8]{ILIE2013}\label{LEMMA8888} Let $Y\in\Omega_0(H)$, and $\alpha:\text{Aff}(Y)\rightarrow G$ be an affine map such that $\alpha|_{Y}$ is proper, then $\alpha$ is proper.
\end{enumerate}
\end{remark}

In the next stage of preparation for our results, we bring some facts about extension of a function $u\in B_p(G_0)$ to a function $B_p(G)$, where $G_0\subseteq G$ is an open subgroup of the locally compact group $G$. We refer the interested reader to \cite[Subsection 3.3]{AHSH2021-II} for detailed description.

Let $G_0\subseteq G$, be any subset, and $u:G_0\rightarrow \mathbb{C}$ be a function. Let $u^\circ$ denote the extension of $u$ to $G$ by setting value zero outside of $G_0$, i.e.
\begin{align*}
u^\circ=\left\{
\begin{array}{ll}
u&\text{on}\; G_0\\
0&\text{o.w.}
\end{array}\right. .
\end{align*}

Following lemma has intricate importance in the sequel.
\begin{lemma}{\cite[Lemma 3]{AHSH2021-II}}\label{LEMMARESTRICTIONMAP}
Let $(\pi,E)\in\text{Rep}_p(G)$. Then the restriction of $\pi$ to the open subgroup $G_0$, which is denoted by $(\pi_{G_0},E)$ belongs to $\text{Rep}_p(G_0)$. Moreover, for each $f\in L_1(G_0)$ and each $g\in L_1(G)$, we have the following relations
\begin{align}\label{w0}
\pi_{G_0}(f)=\pi(f^\circ),\quad\text{and}\quad \pi_{G_0}(g|_{G_0})=\pi(g\chi_{G_0}).
\end{align}
\end{lemma}

Next proposition is the main building block of the extension problem.
\begin{proposition}\cite[Proposition 4]{AHSH2021-II}\label{PROPEXREC}
Let $G$ be a locally compact group and $G_0$ be its open
subgroup, and let $(\pi,E)\in \text{Rep}_p(G)$. Then the following statements hold.
\begin{enumerate}
\item\label{PROPEXREC1}
The map $S_{\pi_{G_0}}:PF_{p,\pi_{G_0}}(G_0)\rightarrow PF_{p,\pi}(G)$ defined via $S_{\pi_{G_0}}(\pi_{G_0}(f))=\pi(f^\circ)$, for $f \in L_1(G_0)$ is an isometric homomorphism. In fact, we have
the following isometric identification
\begin{align*}
PF_{p,\pi_{G_0}}(G_0) =\overline{\{\pi(f)\ :\ f\in L_1(G),\;\text{supp}(f)\subseteq G_0 \}{}}^{\|\cdot\|_{\mathcal{B}(E)}}\subseteq PF_{p,\pi}(G).
\end{align*}
 \item\label{PROPEXREC2}
The linear restriction mapping $R_\pi : B_{p,\pi}\rightarrow B_{p,\pi_{G_0}}$ which is defined for $u\in B_{p,\pi}$, as $R_\pi(u) =u|_{G_0}$ is the dual map of $S_{\pi_{G_0}}$, and is a quotient
map.
\item\label{PROPEXREC2.5}
The extension map $E_\pi: B_{p,\pi_{G_0}}\rightarrow B_{p,\pi}$, defined via $E_\pi(u)=u^\circ$ is an isometric map.
\item\label{PROPEXREC3}
The restriction mapping $R : B_p(G)\rightarrow B_p(G_0)$ is a contraction.
\item\label{PROPEXREC4}
When $(\pi,E)$ is also a $p$-universal representation, we have the following contractive inclusions
\begin{align*}
PF_p(G_0)^*\subseteq B_{p,\pi_{G_0}}\subseteq B_p(G_0)\subseteq\mathcal{M}(A_p(G_0)).
\end{align*}
Under the assumption that $G_0$ is amenable, we have isometric identification below
\begin{align*}
PF_p(G_0)^* = B_{p,\pi_{G_0}}= B_p(G_0)=\mathcal{M}(A_p(G_0)).
\end{align*}
\end{enumerate}
\end{proposition}

Next theorem has a key role in studying of weighted homomorphisms \eqref{eq1}, and it is stated as Proposition 5 in \cite{AHSH2021-II}.

\begin{theorem}[Extension Theorem]\label{PROPEXTENSION}
Let $G$ be a locally compact group and $G_0$ be its open
subgroup. Then
\begin{enumerate}
\item\label{PROPEXTENSION1}
the extension mapping $E_{MM} : \mathcal{M}(A_p(G_0))\rightarrow \mathcal{M}(A_p(G))$, defined for
$u \in \mathcal{M}(A_p(G_0))$ via $E_{MM} (u) = u^\circ$
is an isometric map.
\item\label{PROPEXTENSION2}
for every $u\in B_p(G_0)$, we have $u^\circ\in \mathcal{M}(A_p(G))$, and the map $E_{BM} :
B_p(G_0)\rightarrow \mathcal{M}(A_p(G))$, with $u\mapsto u^\circ$, is a contraction.
\item\label{PROPEXTENSION3}
if $G_0$ is also an amenable subgroup, then for every $u \in B_p(G_0)$, we have $u^\circ \in B_p(G)$, and the associated extending map $E_{BB} : B_p(G_0)\rightarrow B_p(G)$ is an isometric one.
\end{enumerate}
\end{theorem}

\begin{remark}\label{REMARKEXTENSION}
It is worthwhile to take notice of the fact that when the subgroup $G_0\subseteq G$ is an amenable open one, a $p$-universal representation of $G_0$ can be induced by restriction of a $p$-universal representation of $G$ to $G_0$.
\end{remark}
Next and final theorem of this section is a consequence of the Theorem \ref{PROPEXTENSION} and can be found in \cite{AHSH2021-II}.

\begin{theorem}\label{THEOREMCONCLUSION}
Let $G$ and $H$ be locally compact groups, and $\alpha :Y=\cup_{k=1}^nY_k\subseteq H\rightarrow G$ be a continuous piecewise affine map with disjoint $Y_k\in\Omega_{\text{am-}0}(H)$, for $k=1,\ldots,n$. Then  $u\in B_p(G)$ implies that $(u\circ \alpha)^\circ\in B_p(H)$, and consequently, the weighted homomorphism $\Phi_\alpha:B_p(G)\rightarrow B_p(H)$ is well-defined bounded homomorphism.
\end{theorem}

\section{Special $p$-completely bounded operators on $B_p(G)$}\label{SECTIONSPECIALMAPS}

The essential features of the $p$-analog of the Fourier-Stieltjes algebra and the algebra of $p$-pseudofunctions have been provided so far, and this section is allocated to maintain critical tools in dealing with operators on $B_p(G)$. Upcoming theorem is going to reveal the relation between $p$-universal representations of a locally compact group $G$ and the quotient group $G/N$, where $N$ is a closed normal subgroup. Recall the canonical quotient map below
\begin{align*}
q:G\rightarrow G/N, \quad q(x)=xN,\quad x\in G,
\end{align*}
which is a continuous and onto homomorphism, for a closed normal subgroup $N$ of $G$.
\begin{theorem}\label{PROPQUOTIENT}
 Let $N\subseteq G$ be a closed normal subgroup. Then we have the following statements.
\begin{enumerate}
\item\label{PROPQUOTIENT1}
$(\rho,F)\in\text{Rep}_p(G/N)$ implies that $(\rho\circ q,F)\in\text{Rep}_p(G)$ and the identification 
\begin{align*}
PF_{p,\rho}(G/N)=PF_{p,\rho\circ q}(G),
\end{align*}
 is an isometric isomorphism.
 \item\label{PROPQUOTIENT2}
 Each representation $(\pi,E)\in\text{Rep}_p(G)$ induces a representation $(\tilde{\pi}_K,K)\in\text{Rep}_p(G/N)$ so that
 \begin{align*}
 (\tilde{\pi}_K\circ q,K)\subseteq(\pi,E),
 \end{align*}
 where
\begin{align}\label{K-SPACE}
K=\big\{\xi\in E\ :\ \pi(n)\xi=\xi,\; \text{for all}\; n\in N\big\}.
\end{align}
Moreover, if the representation $(\pi,E)$ is $p$-universal, then the induced representation $(\tilde{\pi}_K,K)$ is $p$-universal, as well. 
\item\label{PROPQUOTIENT3}
Denoting by $B_p(G:N)$ the subalgebra of functions in $B_p(G)$ those are constant on each coset of $N$, we have
\begin{align*}
B_p(G:N)=B_p(G/N).
\end{align*}
\item\label{PROPQUOTIENT4}
The map $\Phi_q:B_p(G/N)\rightarrow B_p(G)$ defined through $\Phi_q(u)=u\circ q$, for $u\in B_p(G/N)$ is an isometric map.
\end{enumerate}
\begin{proof}
\begin{enumerate}
\item
Let $(\rho,F)\in\text{Rep}_p(G/N)$. Evidently, we have $(\rho\circ q,F)\in\text{Rep}_p(G)$. Recall the natural map $P:L_1(G)\rightarrow L_1(G/N)$ from \cite{FOLLAND1995},
\begin{align}\label{P-MAP}
Pf(xN)=\int_Nf(xn)dn,\quad xN\in G/N.
\end{align}
For each $x\in G$, $\zeta\in F$, and $\psi\in F^*$ we have $\langle\rho\circ q(x)\zeta,\psi\rangle =\langle\rho(xN)\zeta,\psi\rangle$, and since for each continuous function $u:G/N\rightarrow\mathbb{C}$, and $f\in L_1(G)$, we have $P((u\circ q)\cdot f)=u\cdot Pf$, it follows that
\begin{align*}
\langle\rho\circ q(f)\zeta,\psi\rangle =\langle\rho(Pf)\zeta,\psi\rangle,\ f\in L_1(G),\; \zeta\in F,\;\psi\in F^*.
\end{align*}
Consequently, we have
\begin{align}\label{K1}
\rho\circ q(x)=\rho(xN),\quad\text{and}\quad \rho\circ q(f)=\rho(Pf),\quad x\in G,\; f\in L_1(G).
\end{align}
The relation \eqref{K1} means that $PF_{p,\rho\circ q}(G)=PF_{p,\rho}(G/N)$, isometrically, by the fact that the map $P$ is onto.
\item
Let $(\pi,E)\in\text{Rep}_p(G)$. Recall the closed subspace $K$ of $E$ in \eqref{K-SPACE}. The space $K$ is a $QSL_p$-space and on top of that, it is invariant under $\pi$, i.e.
\begin{align*}
\pi(x)K\subseteq K,\quad\text{and}\quad \pi(f)K\subseteq K,\qquad x\in G, \ f\in L_1(G).
\end{align*}
So, if $(\pi_K,K)$ denotes the representation of $G$ that maps every element $x\in G$ to $\pi_K(x)=\pi(x)|_K$, then
\begin{align*}
(\pi_K,K)\in\text{Rep}_p(G),\quad\text{and}\quad (\pi_K,K)\subseteq (\pi,E).
\end{align*}
Now, put
\begin{align*}
\tilde{\pi}_K:G/N\rightarrow \mathcal{B}(K),\quad \tilde{\pi}_K(xN)=\pi(x),\quad xN\in G/N.
\end{align*}
Then by the definition of $K$, the pair $(\tilde{\pi}_K,K)$ is well-defined and belongs to $\text{Rep}_p(G/N)$. Additionally,
\begin{align*}
(\tilde{\pi}_K\circ q,K) =(\pi_K,K)\subseteq (\pi,E).
\end{align*}
The idea for the $p$-universality part is based on the fact that each element $(\rho,F)$ in $\text{Cyc}_p(G/N)$ can be corresponded to $(\rho\circ q,F)$ in $\text{Cyc}_p(G)$. Indeed, if $(\rho,F)\in\text{Cyc}_p(G/N)$ is associated with the cyclic vector $\xi_0\in F$ then we shall show that $(\rho\circ q,F)$ is a cyclic representation with the same vector $\xi_0$. To this end, assume that $\xi\in F$ is an arbitrary vector. Then due to the fact that $(\rho,F)$ is cyclic, there exists $(g_i)_{i\in\mathbb{I}}\subseteq L_1(G/N)$ such that
\begin{align*}
    \lim_{i\in\mathbb{I}}\|\rho(g_i)\xi_0-\xi\|=0.
\end{align*}
On the other hand, the map $P$ in \eqref{P-MAP} is onto, therefore, for each $i\in\mathbb{I}$, there exists $ f_i\in L_1(G)$ so that $P(f_i)=g_i$. Thus, via \eqref{K1}, we have
\begin{align*}
    \lim_{i\in\mathbb{I}}\|\rho\circ q (f_i)\xi_0-\xi\|=0.
\end{align*}
Hence, $(\rho\circ q,F)$ is a cyclic representation of $G$ associated with $\xi_0$. Now, if $(\pi, E)$ is a $p$-universal representation of $G$, then by the definition, it contains $(\rho\circ q,F)$, and we have (up to an isometry)
\begin{align*}
F\subseteq E,\quad \rho\circ q(f)=\pi(f)|_F\quad f\in L_1(G).
\end{align*}
Moreover, by the definition of $K$, we have $F\subseteq K$, and it is obtained that
\begin{align*}
(\rho\circ q,F)\subseteq (\pi_K,K),\quad\text{and consequently} \quad (\rho,F)\subseteq (\tilde{\pi}_K,K).
\end{align*}
Now, since $(\tilde{\pi}_K,K)$ contains every cyclic representation $(\rho,F)$ of $G/N$, then it is a $p$-universal representation.
\item
Let 
\begin{align*}
u(x)=\langle\pi(x)\xi_0,\eta_0\rangle,\quad x\in G,
\end{align*}
be a function which is not identically zero and is constant on each coset of $N$, where $(\pi,E)\in\text{Rep}_p(G)$, $\xi_0\in E$ and $\eta_0\in E^*$. Assume that
\begin{align*}
F_{\eta_0}=\overline{\{\pi(f)^*\eta_0\ : \ f\in L_1(G)\}{}}^{\|\cdot\|_{E^*}},
\end{align*}
and then define $E_{\eta_0}=E/F_{\eta_0}^{\bot}$. Consequently, we have $E_{\eta_0}=F_{\eta_0}^*$, or equivalently, $E_{\eta_0}^*=F_{\eta_0}$. So, we can assume that the function $u$ is a coefficient function of the representation $(\pi_{\eta_0 }, E_{\eta_0})=(\pi|_{E_{\eta_0}}, E_{\eta_0})$. Now, similar to \eqref{K-SPACE}, define
\begin{align*}
E_{u}=\left\{ \xi\in E_{\eta_0}\ : \ \pi_{\eta_0 }(n)\xi=\xi ,\; \text{for all} \; n\in N\right\}.
\end{align*}
One can easily check that $\xi_0\in E_u$. Indeed, through the fact that $u\in B_p(G:N)$, for an arbitrary element $f\in L_1(G)$, we have
\begin{align*}
    \langle \pi_{\eta_0 }(n)\xi_0,\pi_{\eta_0 }(f)^*\eta_0\rangle&= \langle \pi_{\eta_0 }(f)\pi_{\eta_0 }(n)\xi_0,\eta_0\rangle\\
    &=\int f(x)\langle \pi_{\eta_0 }(x)\pi_{\eta_0 }(n)\xi_0,\eta_0\rangle dx\\
    &=\int f(x)\langle\pi_{\eta_0 }(xn)\xi_0,\eta_0\rangle dx\\
    &=\int f(x)u(xn) dx\\
    &=\int f(x)u(x) dx\\
    &=\langle  \xi_0,\pi_{\eta_0 }(f)^*\eta_0\rangle .
\end{align*}
This implies that $\pi_{\eta_0}(n)\xi_0=\xi_0$, for all $n\in N$. Following the notations of Part \eqref{PROPQUOTIENT2}, if we put $(\pi_u,E_u)=(\pi_{\eta_0}|_{E_u},E_u)$ then $(\pi_u,E_u)\in\text{Rep}_p(G)$, and the representation $(\tilde{\pi}_u,E_u)$ is well-defined, and belongs to $\text{Rep}_p(G/N)$. On top of that, if we consider
\begin{align}\label{U-TILDE}
\tilde{u}(xN)=\langle\tilde{\pi}_u(xN)\xi_0,\eta_0\rangle,\quad xN\in G/N.
\end{align}
then we have $\tilde{u}\in B_p(G/N)$ and that $u=\tilde{u}\circ q$. Let us define
\begin{align}
& k_N:B_p(G:N)\rightarrow B_p(G/N),\quad k_N(u)=\tilde{u},\label{MAPSKN}\\
& l_N: B_p(G/N)\rightarrow B_p(G:N),\quad l_N(u)=u\circ q.\label{MAPSLN}
\end{align}
where $\tilde{u}$ is as \eqref{U-TILDE}. Evidently, these maps are contractive homomorphisms and are inverse of each other. Therefore, they are isometric isomorphisms.
\item
It is straightforward through previous parts. Indeed, Part \eqref{PROPQUOTIENT2} implies that the map $\Phi_q$ is a contraction, and by Part \eqref{PROPQUOTIENT3} we have the isometry.
\end{enumerate}
\end{proof}
\end{theorem}

Next theorem is our first main result of this paper, and it will be applied to give the results on weighted homomorphisms on the $p$-analog of the Fourier-Stieltjes algebras. For more clarification, we need to introduce the notion of the $p$-tensor product $E\tilde{\otimes}_p F$ of two $QSL_p$-spaces $E$ and $F$, that is defined in \cite{RUNDE2005}. In fact, Runde introduced the norm $\|\cdot\|_p$ on the algebraic tensor product $E\otimes F$ which benefits from pivotal properties. An important property of the norm $\|\cdot\|_p$ is the fact that the completion $E\tilde{\otimes}_p F$ of $E\otimes F$ with respect to $\|\cdot\|_p$ is a $QSL_p$-space. Furthermore, for two representations $(\pi,E)$ and $(\rho,F)$ of the locally compact group $G$ in $\text{Rep}_p(G)$, the representation $(\pi\otimes\rho,E\tilde{\otimes}_p F)$ is well-defined and belongs to $\text{Rep}_p(G)$. As a result, for two functions $u(x)=\langle\pi(x)\xi,\eta\rangle$ and $v(x)=\langle\rho(x)\zeta,\psi\rangle,\quad (x\in G)$, the pointwise product of them is a coefficient function of the representation $(\pi\otimes\rho, E\tilde{\otimes}_p F)$, i.e.,
\begin{align*}
u\cdot v(x)=\langle(\pi(x)\otimes\rho(x))(\xi\otimes\zeta),\eta\otimes\psi\rangle,\quad x\in G .    
\end{align*}
For more details on $p$-tensor product $\tilde{\otimes}_p$ see \cite[Theorem 3.1 and Corollary 3.2]{RUNDE2005}.

\begin{theorem}\label{IMPORTANTMAPPING} Let $p\in (1,\infty)$ and $G$ be a locally compact group. Then we have the following statements.
\begin{enumerate}
\item\label{ITEMIDENTITYMAP}
For any $(\rho, F)\in \text{Rep}_{p}(G)$, the identity map $I:B_{p,\rho}\rightarrow B_p(G)$ is a $p$-completely contractive map.
\item\label{ITEMRESTRICTIONMAP}
For an open subgroup $G_0$ of $G$, the restriction map $R: B_p(G)\rightarrow B_p(G_0)$, is a $p$-completely contractive homomorphism.
\item\label{ITEMTRANSLATIONMAP}
For an element $a\in G$, the translation map $L_a:B_p(G)\rightarrow B_p(G)$, defined through $L_a(u)={}_a u$, where ${}_au(x)=u(ax)$, for $x\in G$, is a $p$-completely contractive map.
\item\label{ITEMQUOTIENTMAP}
The homomorphism $\Phi_q: B_p(G/N)\rightarrow B_p(G)$, with $\Phi_q(u)= u\circ q$, is a $p$-completely contractive homomorphism.
\item\label{ITEMEXTENSIONMAP}
For an open amenable subgroup $G_0$ of $G$, the extension map $E_{BB}:B_p(G_0)\rightarrow B_p(G)$ is a $p$-completely contractive homomorphism.
\item\label{ITEMMULTIPLICATIONYYY}
For an open coset $Y$ of an open subgroup $G_0$ of $G$, the map $M_Y:B_p(G)\rightarrow B_p(G)$, with $M_Y(u)=u\cdot\chi_Y$, is $p$-completely contractive homomorphism. More generally, for a set $Y\in\Omega_0(G)$, the map $M_Y$ is a $p$-completely bounded homomorphism.
\end{enumerate}

\begin{proof}

\begin{enumerate}

\item
We want to prove that for each $(\rho,F)\in\text{Rep}_{p}(G)$, the following map is a $p$-complete contraction,
\begin{align}\label{ID}
I: B_{p,\rho}\rightarrow B_p(G),\quad I(u)=u .
\end{align}

Let $(\pi ,E)$ be a $p$-universal representation of $G$ that contains the representation $(\rho,F)$. Following relations hold between $(\rho ,F)$, and $(\pi, E)$,
\begin{align*}
F\subseteq E, \quad \rho(x)=\pi(x)|_{F},\quad\text{and}\quad\rho(f)=\pi(f)|_{F},\quad x\in G,\; f\in L_1(G).
\end{align*}
Since $\rho(f)=\pi(f)|_{F}$, then $\|\rho(f)\|\leq\|\pi(f)\|$. Additionally, the map $I$ is weak$^*$-weak$^*$ continuous, and it is a contraction by Theorem \ref{THEOREMRUNDE2005}.
Define
\begin{align*}
{}_*I: UPF_{p}(G)\rightarrow PF_{p,\rho}(G),\quad {}_*I(\pi(f))=\pi(f)|_{F}=\rho(f),
\end{align*}
then ${}_*I$ is the predual of the map \eqref{ID}, since for every $f\in L_1(G)$ and $u\in B_{p,\rho}$ we have $\langle\pi(f), I(u)\rangle=\langle\rho(f),u\rangle$. The following calculations indicate that ${}_*I$ is a $p$-complete contraction: for each $n\in\mathbb{N}$, and $[\rho(f_{ij})]\in\mathbb{M}_n(PF_{p,\rho}(G))$ we have
\begin{align*}
\|[\rho(f_{ij})]\|_n&=\sup\bigg\{\|[\rho(f_{ij})](\xi_{j})_{j=1}^n\|\ :\  (\xi_j)_{j=1}^n\subseteq F,\; \sum_{j=1}^n\|\xi_j\|^p\leq 1\bigg\}\\
&=\sup\bigg\{\|[\pi(f_{ij})](\xi_{j})_{j=1}^n\|\ :\  (\xi_j)_{j=1}^n\subseteq F,\; \sum_{j=1}^n\|\xi_j\|^p\leq 1\bigg\}\\
&\leq\sup\bigg\{\|[\pi(f_{ij})](\xi_{j})_{j=1}^n\|\ :\  (\xi_j)_{j=1}^n\subseteq E,\; \sum_{j=1}^n\|\xi_j\|^p\leq 1\bigg\}\\
&=\|[\pi(f_{ij})]\|_n,
\end{align*}
so, we have $\|[\rho(f_{ij})]\|_n\leq \|[\pi(f_{ij})]\|_n$, and by this, it is concluded that
\begin{align*}
\|I\|_{\text{p-cb}}\leq \|{}_*I\|_{\text{p-cb}}\leq 1 .
\end{align*} 

\item
Let $G_0\subseteq G$ be an open subgroup and recall the maps $R$ and $S=S_\pi$ (for a $p$-universal representation $(\pi,E)$ of $G$) in Proposition \ref{PROPEXREC} together with the relation \eqref{w0} from Lemma \ref{LEMMARESTRICTIONMAP}. Since we have $S^*=R$, and the map $S$ is indeed an identity, therefore, $S$ is a $p$-completely isometric map, and consequently, 
\begin{align*}
\|R\|_{\text{p-cb}}=\|S^*\|_{\text{p-cb}}\leq\|S\|_{\text{p-cb}}=1.
\end{align*}

\item
For $a\in G$, consider the following map
\begin{align*}
L_a:B_p(G)\rightarrow B_p(G),\quad L_a(u)={}_au,\quad {}_au(x)=u(ax),\; x\in G.
\end{align*}
The predual of the map of $L_a$ is as following
\begin{align*}
&{}_*L_a :UPF_{p}(G)\rightarrow UPF_p(G),\quad {}_*L_a(\pi(f))=\pi(\lambda_p(a)f)),
\end{align*}
and it is clearly $p$-completely contractive, consequently, this holds true for $L_a$. On the other hand, the map $L_a$ has the inverse $L_{a^{-1}}$, and it is $p$-completely contractive as well, which makes $L_a$ a $p$-completely isometric map.

\item Let $N\subseteq G$ be a closed normal subgroup. Recalling all notations from Theorem \ref{PROPQUOTIENT}, we let $(\pi,E)\in\text{Rep}_p(G)$ be a $p$-universal representation and $(\tilde{\pi}_K,K)$ be the induced $p$-universal representation of $G/N$. 
For functions $f\in L_1(G)$, and $u\in B_p(G/N)$, we have
\begin{align}\label{EQPHQ}
\langle\pi(f),\Phi_q(u)\rangle=\langle\pi(f),u\circ q\rangle=\langle\tilde{\pi}(Pf),u\rangle.
\end{align}
This implies that the map $\Phi_q$ is weak$^*$-weak$^*$ continuous, and by this we define the predual map ${}_*\Phi_q$, as following,
\begin{align*}
{}_*\Phi_q : PF_{p,\tilde{\pi}\circ q}(G)\rightarrow UPF_{p}(G/N),\quad {}_*\Phi_q(\tilde{\pi}\circ q(f))=\tilde{\pi}(P f),\quad f\in L_1(G),
\end{align*}
which by \eqref{EQPHQ} we have $({}_*\Phi_q)^*=\Phi_q$. By using similar relation to \eqref{K1}, we have $\tilde{\pi}\circ q(f)=\tilde{\pi}(Pf)$, which means that the predual map ${}_*\Phi_q$ is an identity that is $p$-completely isometric map via the following computation
\begin{align*}
\|{}_*\Phi_q^{(n)}([\tilde{\pi}\circ q(f_{ij})])\|_n=
\|[\tilde{\pi}(Pf_{ij})]\|_n=
\|[\tilde{\pi}\circ q(f_{ij})]\|_n.
\end{align*}
Therefore, we have $\|\Phi_q\|_{\text{p-cb}}\leq 1$.

\item 
Let $G_0\subseteq G$, be an open amenable subgroup, and $u\in B_p(G_0)$. By Theorem \ref{PROPEXTENSION}, the extension map $E_{BB}$ is well-defined,
\begin{align*}
&E_{BB}: B_p(G_0)\rightarrow B_p(G),\quad E_{BB}(u)=u^\circ .
\end{align*}
Let $(\pi ,E)$ be a $p$-universal representation of $G$. We denote the restriction of $(\pi,E)$ to $G_0$ by $(\pi_{G_0}, E)$, which is a $p$-universal representation of $G_0$ via Remark \ref{REMARKEXTENSION}.
We note that by the relation
\begin{align}\label{EQEXT}
\langle \pi(f), u^\circ\rangle=\langle\pi_{G_0}(f|_{G_0}),u\rangle,\quad f\in L_1(G),\ u\in B_p(G_0),
\end{align}
the map $E_{BB}$ is weak$^*$-weak$^*$ continuous. So, we define the predual map ${}_*E_{BB}$, as following
\begin{align*}
{}_*E_{BB}: UPF_p(G)\rightarrow UPF_p(G_0),\quad {}_*E_{BB}(\pi(f)):=\pi_{G_0}(f|_{G_0}),
\end{align*}
which by \eqref{EQEXT} we have $({}_*E_{BB})^*=E_{BB}$.
We need to take notice of the fact that since $\chi_{G_0}\in B_p(G)$, via \cite[Theorem 2]{AHSH2021-II}, $\chi_{G_0}$ is a normalized coefficient function of $(\pi, E)$, i.e. there exist $\xi_\chi\in E$, and $\eta_\chi\in E^*$ so that
\begin{align}\label{p0}
\|\xi_\chi\|=\|\eta_\chi\|=1, \ \text{and}\quad \chi_{G_0}(x)=\langle\pi(x)\xi_\chi,\eta_\chi\rangle,\quad x\in G, 
\end{align}
On the other hand, for $f\in L_1(G_0)$, $\xi\in E$, and $\eta\in E^*$, we have
\begin{align*}
\langle\pi(f\chi_{G_0})\xi,\eta\rangle&=\int_{G}f(x)\chi_{G_0}(x)\langle\pi(x)\xi,\eta\rangle dx\\
&=\int_{G}f(x)\langle\pi(x)\xi_\chi,\eta_\chi\rangle\langle\pi(x)\xi,\eta\rangle dx\\
&=\int_{G}f(x)\langle(\pi(x)\otimes\pi(x))(\xi_\chi\otimes\xi),\eta_\chi\otimes\eta\rangle dx\\
&=\langle(\pi\otimes\pi(f))(\xi_\chi\otimes\xi),\eta_\chi\otimes\eta\rangle,
\end{align*}
which implies that
\begin{align}\label{Z1}
\langle\pi(f\chi_{G_0})\xi,\eta\rangle=\langle(\pi\otimes\pi(f))(\xi_\chi\otimes\xi),\eta_\chi\otimes\eta\rangle,\quad f\in L_1(G),\ \xi\in E,\ \eta\in E^*.
\end{align}
Therefore, by combining equality \eqref{Z1} with \eqref{w0}, we have
\begin{align}\label{EQNEW}
\langle\pi_{G_0}(f|_{G_0})\xi,\eta\rangle=\langle(\pi\otimes\pi(f))(\xi_\chi\otimes\xi),\eta_\chi\otimes\eta\rangle,\quad f\in L_1(G),\ \xi\in E,\ \eta\in E^*.
\end{align}
Additionally, since $(\pi,E)$ is a $p$-universal representation, and we have
\begin{align*}
(\pi, E)\subseteq (\pi\otimes\pi, E\tilde{\otimes}_p E),
\end{align*}
thus, $(\pi\otimes\pi, E\tilde{\otimes}_p E)$ can be assumed as a $p$-universal of $G$. Let
\begin{align*}
{}_*E_{BB}^{(n)}: \mathbb{M}_n(UPF_p(G))\rightarrow \mathbb{M}_n(UPF_p(G_0)),\quad {}_*E_{BB}^{(n)}(\pi(f_{ij})):=(\pi_{G_0}(f_{ij}|_{G_0})).
\end{align*}
Then via \eqref{EQNEW} we have

\begin{small}

\begin{align*}
\|{}_*E_{BB}^{(n)}([\pi(f_{ij})])\|_n&=\|[\pi_{G_0}(f_{ij}|_{G_0})]\|_n\\
&=\sup\bigg\{|\sum_{i,j=1}^n\langle\pi_{G_0}(f_{ij}|_{G_0})\xi_j,\eta_i\rangle|\ :\ (\xi_j)_{j=1}^n\subseteq E, \\
&\qquad\qquad (\eta_i)_{i=1}^n\subseteq E^*,\;\sum_{j=1}^n\|\xi_j\|^p\leq 1,\;  \sum_{i=1}^n\|\eta_i\|^{p'}\leq 1\bigg\}\\
&\leq\sup\bigg\{|\sum_{i,j=1}^n\langle(\pi\otimes\pi(f_{ij}))\phi_j,\psi_i\rangle|\ :\ (\phi_j)_{j=1}^n\subseteq E\tilde{\otimes}_pE,\\
&\qquad\qquad(\psi_i)_{i=1}^n\subseteq E^*\tilde{\otimes}_{p'}E^*,\;\sum_{j=1}^n\|\phi_j\|^p\leq 1,\;  \sum_{i=1}^n\|\psi_i\|^{p'}\leq 1\bigg\}\\
&=\|(\pi\otimes\pi(f_{ij}))\|_n,
\end{align*}
\end{small}
and since norm of $UPF_p(G)$ is independent of choosing $p$-universal representation (see Theorem \ref{TH45}) then we have $\|{}_*E_{BB}\|_{\text{p-cb}}\leq 1$, which implies that $\|E_{BB}\|_{\text{p-cb}}\leq 1$.

\item
Let $(\pi,E)$ be a $p$-universal representation. By \cite[Corollary 2]{AHSH2021-II}, for $Y\in\Omega_0(G)$ the map $M_Y: B_p(G)\rightarrow B_p(G)$ with $M_Y(u)=u\cdot \chi_Y$ is well-defined, and 
\begin{align*}
\|M_Y\|\leq 2^{m_Y},
\end{align*}
On the other hand, by the following relation this map is weak$^*$-weak$^*$ continuous
\begin{align}\label{EQMULTBP}
\langle\pi(f),u\cdot\chi_Y\rangle=\langle\pi(f\cdot\chi_Y),u\rangle,\quad f\in L_1(G),\ u\in B_p(G).
\end{align}
So, one may define its predual map as following
\begin{align*}
{}_*M_Y:UPF_p(G)\rightarrow UPF_p(G),\quad {}_*M_Y(\pi(f))=\pi(f\cdot\chi_Y),
\end{align*} 
and by \eqref{EQMULTBP} we have $({}_*M_Y)^*=M_Y$. 
\begin{enumerate}

\item[Step 1:] To prove the claim, first we let $Y$ be an open coset itself. Similar to \eqref{p0}, the function $\chi_Y$ is a normalized coefficient function of the representation $(\pi,E)$ which means  that there are elements $\xi_Y\in E$, and $\eta_Y\in E^*$ such that
\begin{align*}
\|\xi_Y\|=\|\eta_Y\|=1,\ \text{and}\quad \chi_Y(x)=\langle\pi(x)\xi_Y,\eta_Y\rangle,\quad x\in G.
\end{align*}
So, for a matrix $[\pi(f_{ij})]\in \mathbb{M}_n(UPF_p(G))$, through the relation \eqref{Z1}, we have
\begin{small}

\begin{align*}
\|[\pi(f_{ij}\cdot\chi_Y)]\|_n&=\sup\bigg\{|\sum_{i,j=1}^n\langle\pi(f_{ij}\cdot\chi_Y)\xi_j,\eta_i\rangle|\ :\ (\xi_j)_{j=1}^n\subseteq E,\\
&\qquad\qquad (\eta_i)_{i=1}^n\subseteq E^*,\; \sum_{j=1}^n\|\xi_j\|^p\leq 1,\; \sum_{i=1}^n\|\eta_i\|^{p'}\leq 1 \bigg\}\\
&\leq \sup\bigg\{|\sum_{i,j=1}^n\langle\pi\otimes\pi(f_{ij})\phi_j,\psi_i\rangle|\ :\ (\phi_j)_{j=1}^n\subseteq E\tilde{\otimes}_pE,\\
&\qquad\qquad (\psi_i)_{i=1}^n\subseteq E^*\tilde{\otimes}_{p'}E^*,\; \sum_{j=1}^n\|\phi_j\|^p\leq 1,\; \sum_{i=1}^n\|\psi_i\|^{p'}\leq 1 \bigg\}\\
&=\|(\pi\otimes\pi(f_{ij}))\|_n.
\end{align*}

\end{small}
By these computations, we obtain that the map ${}_*M_Y$ is a $p$-complete contraction. Therefore, we have $\|M_Y\|_{\text{p-cb}}\leq 1$, through the fact that $p$-operator norm of $UPF_{p}(G)$ is independent of choosing $p$-universal representation.

\item[Step 2:] Now let $Y=Y_0\backslash\cup_{i=1}^m Y_i \in\Omega_0(G)$, and we have,
\begin{align*}
M_{Y}=M_{{Y_0}}-(\sum_{i=1}^m M_{Y_i} -\sum_{i,j}M_{{Y_i\cap Y_j}}+\sum_{i,j,k}M_{{Y_i\cap Y_j\cap Y_k}}+\ldots+(-1)^{m+1}M_{{Y_1\cap\ldots Y_m}}).
\end{align*}
Therefore, we have $\|M_Y\|_{\text{p-cb}}\leq 2^{m_Y}$.

\end{enumerate}
\end{enumerate}
\end{proof}
\end{theorem}
The following, is an immediate consequence of the Theorem \ref{PROPQUOTIENT} and Theorem \ref{IMPORTANTMAPPING}-\eqref{ITEMQUOTIENTMAP}.
\begin{corollary}
For a closed normal subgroup $N$ of a locally compact group $G$ the identification $B_p(G/N)=B_p(G:N)$ is a $p$-complete isometry.
\begin{proof}
Recall the maps $k_N$ and $l_N$ in \eqref{MAPSKN} and \eqref{MAPSLN}. Evidently, they are weak$^*$-weak$^*$ continuous. Then, with notations of Theorem \eqref{PROPQUOTIENT}, let $(\pi,E)$ be a $p$-universal representation of $G$. It should be noted that elements of the algebra $B_p(G:N)$ are coefficient functions of the representation $(\tilde{\pi}\circ q, K)$. Hence one my define their predual maps as following
\begin{align*}
    &{}_*k_N: UPF_p(G/N)\to PF_{p,\tilde{\pi}\circ q}(G),\quad {}_*k_N(\tilde{\pi}(f))=\tilde{\pi}\circ q(g),\quad f\in L_1(G/N),\\
    &\text{where}\; g\in L_1(G),\; \text{with}\; P(g)=f,\\
    &{}_*l_N: PF_{p,\tilde{\pi}\circ q}(G) \to UPF_p(G/N),\quad {}_*l_N(\tilde{\pi}\circ q (g))=\tilde{\pi} (P(g)),\quad g\in L_1(G ).
\end{align*}
By \eqref{K1} the map ${}_*k_N$ is well-defined, and in fact, both maps ${}_*k_N$ and ${}_*l_N$ are identities and inverses of each other, from which the result follows.
\end{proof}
\end{corollary}

\begin{remark}
\begin{enumerate}

The importance of Theorem \ref{IMPORTANTMAPPING}-\eqref{ITEMIDENTITYMAP} is that while we are working with maps with ranges as subspaces of the $p$-analog of the Fourier-Stieltjes algebras, we just need to restrict ourselves to their ranges, as what we have done in the rest of Theorem \ref{IMPORTANTMAPPING}.

\end{enumerate}
\end{remark}

\section{$p$-Completely bounded homomorphisms on $B_p(G)$ induced by proper  piecewise affine maps}\label{SECTIONINDUCEDHOMO}

As an application of previous sections, we are ready to study on homomorphisms $\Phi_\alpha : B_p(G)\rightarrow B_p(H)$ of the form
\begin{align*}
\Phi_\alpha u= \left\{
\begin{array}{ll}
u\circ \alpha & \text{on} \: Y\\
0 & \text{o.w}
\end{array}\right.,\quad u\in B_p(G),
\end{align*}
for the proper  and continuous piecewise affine map $\alpha :Y\subseteq H\rightarrow G$ with $Y=\cup_{i=1}^nY_i$ and $Y_i\in\Omega_{\text{am-}0}(H)$, which are pairwise disjoint, for $i=1,\ldots,n$. We will give some results in the sequel, and we need the following lemma in this regard. For general form of this lemma, see \cite[Lemma 1]{ILIE2013}, and related references there, e.g. \cite{DERIGHETTI1982}.
\begin{lemma}\label{LEMMACALPHA}
Let $G$ and $H$ be locally compact groups and $\alpha :H\rightarrow G$ be a proper homomorphism that is onto, then there is a constant $c_\alpha>0$, such that
\begin{align*}
\int_H f\circ\alpha(h)dh=c_\alpha\int_G f(x)dx,\quad f\in L_1(G).
\end{align*}

\end{lemma}

\begin{proposition}\label{ITEMHOMOMORPHISMPHI}
Let $G$ and $H$ be locally compact groups and $\alpha :H\rightarrow G$ be a proper continuous group homomorphism. Then the homomorphism $\Phi_\alpha : B_p(G)\rightarrow B_p(H)$, of the form $\Phi_\alpha(u)=u\circ\alpha$, is well-defined and $p$-completely contractive homomorphism.

\begin{proof}

Let $(\pi,E)$ be a $p$-universal representation of $G$. Obviously, $(\pi\circ\alpha , E)\in\text{Rep}_p(H)$, and $\Phi_\alpha$ is a contractive map so that its range is the subspace of $B_p(H)$ of functions which are coefficient functions of the representation $(\pi\circ\alpha , E)$, i.e. the space $B_{p,\pi\circ\alpha}$. We divide our proof into two steps.

\begin{enumerate}
\item[Step 1:]
First, we suppose that $\alpha :H\rightarrow G$ is a continuous isomorphism.
In this case, $(\pi\circ\alpha , E)$ is a $p$-universal representation of $H$, and by Lemma \ref{LEMMACALPHA},  for every $f\in L_1(H)$ and $u\in B_p(G)$, we have
\begin{align*}
\langle\pi\circ\alpha(f),u\circ\alpha\rangle&=\int_H f(h)u\circ\alpha(h)dh\\
&=\int_H (f\circ\alpha^{-1})\circ\alpha(h)u\circ\alpha(h)dh\\
&=c_\alpha\int_G f\circ\alpha^{-1}(x)u(x)dx\\
&=c_\alpha\langle\pi(f\circ\alpha^{-1}),u\rangle.
\end{align*}
Consequently, the map $\Phi_\alpha$ is weak$^*$-weak$^*$ continuous, and we define
\begin{align*}
{}_*\Phi_\alpha :UPF_p(H)\rightarrow UPF_p(G),\quad{}_*\Phi_\alpha(\pi\circ\alpha(f)):=c_\alpha\pi(f\circ\alpha^{-1}).
\end{align*}
According to the above relation, we have $({}_*\Phi_\alpha)^*=\Phi_\alpha$. On the other hand, for every $\xi\in E$ and $\eta\in E^*$, we have
\begin{align*}
\langle\pi\circ\alpha(f)\xi,\eta\rangle & =\int_{H}f(h)\langle\pi\circ\alpha(h)\xi,\eta\rangle dh\\
& =\int_{H}f\circ\alpha^{-1}\circ\alpha(h)\langle\pi\circ\alpha(h)\xi,\eta\rangle dh\\
& =c_\alpha\int_{G}f\circ\alpha^{-1}(x)\langle\pi(x)\xi,\eta\rangle dx\\
&=\langle c_\alpha\pi(f\circ\alpha^{-1})\xi,\eta\rangle,
\end{align*}
which means $\pi\circ\alpha (f)=c_\alpha\pi(f\circ\alpha^{-1})$.
Consequently, ${}_*\Phi_\alpha$ is an identity map, so is a $p$-complete isometric map,
\begin{align*}
\|{}_*\Phi_\alpha^{(n)}(\pi\circ\alpha(f_{ij}))\|_n & =\|(c_\alpha\pi(f_{ij}\circ\alpha^{-1}))\|_n =\|(\pi\circ\alpha(f_{ij}))\|_n.
\end{align*}
Therefore, $\|\Phi_\alpha\|_{\text{p-cb}}\leq\|{}_*\Phi_\alpha\|_{\text{p-cb}}=1$.

\item[Step 2:]
Now let $\alpha :H\rightarrow G$ be any proper continuous homomorphism. Let $G_0=\alpha(H)$, and $N=\ker\alpha$. Let us define
\begin{align*}
\tilde{\alpha} :H/N\rightarrow G_0,\quad \tilde{\alpha}(xN)=\alpha(x),\quad x\in G,
\end{align*}
then by Proposition \ref{ramirez}, the map $\tilde{\alpha}$ is a continuous isomorphism, $N$ is a compact normal subgroup of $H$, and $G_0$ is an open subgroup of $G$. Therefore, $\alpha=\tilde{\alpha}\circ q $. By Step 1,  the map $\Phi_{\tilde{\alpha}}$ is $p$-completely contractive, and because of the following composition, $\Phi_\alpha$ is $p$-completely contractive, via Theorem \ref{IMPORTANTMAPPING}-\eqref{ITEMRESTRICTIONMAP}-\eqref{ITEMQUOTIENTMAP},
\begin{align*}
\Phi_\alpha=\Phi_q\circ\Phi_{\tilde{\alpha}}\circ R .
\end{align*}

\end{enumerate}

\end{proof}

\end{proposition}
For the next theorem, we have to put the amenability assumption on the subgroups of $H$, because of Theorem \ref{PROPEXTENSION}.

\begin{theorem}\label{ITEMAFFINE}
Let $G$ and $H$ be two locally compact groups, $Y$ be an open coset of an open amenable subgroup of $H$, and $\alpha :Y\subseteq H\rightarrow G$ be a continuous proper affine map. Then the map $\Phi_\alpha : B_p(G)\rightarrow B_p(H)$, defined as
\begin{align}\label{Z2}
\Phi_\alpha(u)=\left\{
\begin{array}{ll}
u\circ\alpha, & \text{on}\; Y,\\
0,&\text{o.w.}
\end{array}\right.,\quad u\in B_p(G),
\end{align}
is $p$-completely contractive. More generally, if $\alpha$ is a continuous proper piecewise affine map, and $Y=\cup_{i=1}^nY_i$, where disjoint sets $Y_i$ belong to $\Omega_{\text{am-}0}(H)$, then the map $\Phi_\alpha$ is $p$-completely bounded.

\begin{proof}
Through Theorem \ref{THEOREMCONCLUSION} the map \eqref{Z2} is well-defined and it is contractive when the map $\alpha$ is affine on an open coset $Y$ of an open amenable subgroup $H_0$ of $H$, and is bounded in the case that the map $\alpha$ is piecewise affine on the set $Y=\cup_{i=1}^nY_i$, with disjoint $Y_i\in\Omega_{\text{am-}0}(H)$, for $i=1,\ldots, n$, of the form $Y_i=Y_{i,0}\backslash\cup_{j=1}Y_{ij}^{m_i}$.

Let $\alpha :Y=y_0H_0\rightarrow G$ be a continuous proper affine map on the open coset $Y=y_0H_0$, and $H_0$ be an open amenable subgroup of $H$, for which by Remark \ref{AFFFIINEREMMM}-\eqref{affine-remark}, there exists a continuous group homomorphism $\beta :H_0\subseteq H\rightarrow G$ associated with $\alpha$ such that
\begin{align*}
\beta(h)=\alpha(y_0)^{-1}\alpha(y_0h),\quad h\in H_0,
\end{align*}
and it is proper via Remark \ref{AFFFIINEREMMM}-\eqref{HOMAFFPROPER}. Now, consider the following composition
\begin{align*}
\Phi_\alpha= L_{{y_0}^{-1}}\circ E_{BB}\circ\Phi_\beta\circ L_{\alpha(y_0)},
\end{align*}
 then by Proposition \ref{ITEMHOMOMORPHISMPHI}, and Theorem \ref{IMPORTANTMAPPING}-\eqref{ITEMTRANSLATIONMAP}-\eqref{ITEMEXTENSIONMAP}, in place, the map $\Phi_\alpha$ is $p$-completely contractive homomorphism.

Next, we consider the piecewise affine case. Let the map $\alpha :Y\subseteq H\rightarrow G$ be a proper and continuous piecewise affine map. Then for some $n\in\mathbb{N}$, and $i=1,\ldots,n$, there are disjoint sets $Y_i\in\Omega_{\text{am-}0}(H)$, such that $Y=\cup_{i=1}^{n}Y_i$, and  for $i=1,\ldots,n$, the map $\alpha_i:{Aff(Y_i)}\rightarrow G$ is an affine map, together with the fact that $\alpha_i|_{Y_i}=\alpha|_{Y_i}$. On top of that, by Remark \ref{AFFFIINEREMMM}-\eqref{LEMMA8888}, each affine map $\alpha_i$ is proper. Therefore, by considering
\begin{align*}
\Phi_\alpha= \sum_{i=1}^n M_{Y_i}\circ\Phi_{\alpha_i},
\end{align*}
and through the above computations for the maps $\Phi_{\alpha_i}$, we have
\begin{align*}
\|\Phi_\alpha\|_{\text{p-cb}}\leq \sum_{i=1}^n 2^{m_{Y_i}}.
\end{align*}
Here the value $m_{Y_i}$ is the corresponding number to each $Y_i$, as it is in Theorem \ref{IMPORTANTMAPPING}-\eqref{ITEMMULTIPLICATIONYYY}.
\end{proof}
\end{theorem}


\end{document}